\documentclass[12pt]{article}
\renewcommand{\emph}[1]{\textit{#1}}
\usepackage{enumerate,amsmath,latexsym,amssymb}
\usepackage{color}\usepackage{graphicx}

\definecolor{brown}{cmyk}{0, 0.72, 1, 0.45}
\definecolor{grey}{gray}{0.5}

\newcounter{rot}

\newcommand{\ignore}[1]{}

\def\cA{{\cal A}}

\newcommand{\set}[1]{\left\{#1\right\}}

\newcommand{\proofstart}{{\noindent \bf Proof\hspace{2em}}}
\newcommand{\proofend}{\hspace*{\fill}\mbox{$\Box$}\\ \medskip\\ \medskip}

\def\ii_(#1,#2){i_{#1}^{#2}}

\def\LL{\Lambda}

\def\a{\alpha}
\def\b{\beta}

\def\e{\varepsilon}
\def\f{\phi}

\def\l{\lambda}

\def\p{\pi}

\def\s{\sigma}
\def\t{\tau}
\def\om{\omega}

\def\cE{\mathcal{E}}

\def\cT{{\cal T}}

\def\Re{\mathbb{R}}
\newcommand{\brac}[1]{\left( #1 \right)}

\newcommand{\expect}{\operatorname{\bf E}}
\def\E{\expect}
\renewcommand{\Pr}{\operatorname{\bf Pr}}
\newcommand\bfrac[2]{\left(\frac{#1}{#2}\right)}

\newtheorem{theorem}{Theorem}
\newtheorem{lemma}{Lemma}

\newtheorem{remthm}[lemma]{Remark}

\newcounter{thmtemp}

\newcommand{\nospace}[1]{}

\def\path{\operatorname{PATH}}

\newcommand{\beq}[1]{\begin{equation}\label{#1}}
\def\eeq{\end{equation}}

\renewcommand{\Re}{\mathbb{R}}

\def\leb{\leq_b}
\begin{document}
\title{Long paths in random Apollonian networks}

\author{Colin Cooper\thanks{Department of  Informatics,
King's College, University of London, London WC2R 2LS, UK.
Supported in part by EPSRC grant EP/J006300/1}\and Alan
Frieze\thanks{Department of Mathematical Sciences, Carnegie Mellon
University, Pittsburgh PA15213, USA. Supported in part by NSF grant CCF0502793.}}

\maketitle
\begin{abstract}
We consider the length $L(n)$ of the longest path in a randomly generated Apollonian Network (ApN) $\cA_n$. We show that w.h.p. $L(n)\leq ne^{-\log^cn}$ for any constant $c<2/3$.
\end{abstract}

\section{Introduction}
This paper is concerned with the length of the longest path in a random Apollonian Network (ApN) $\cA_n$. We start with a triangle $T_0=xyz$ in the plane. We then place a point $v_1$ in the centre of this triangle creating 3 triangular faces. We choose one of these faces at random and place a point $v_2$ in its middle. There are now 5 triangular faces. We choose one at random and place a point $v_3$ in its centre. In general, after we have added $v_1,v_2,\ldots,v_1$ there will $2n+1$ triangular faces. We choose one at random and place $v_n$ inside it. The random graph $\cA_n$ is the graph induced by this embedding. It has $n+3$ vertices and $3n+6$ edges.

This graph has been the object of study recently. Frieze and Tsourakakis \cite{FT} studied it in the context of scale free graphs. They determined properties of its degree sequence, properties of the spectra of its adjacency matrix, and its diameter.
Cooper and Frieze \cite{CF}, Ebrahimzadeh, Farczadi, Gao,
Mehrabian, Sato,  Wormald and Zung \cite{Wor} improved the diameter result and determine the diameter asymptotically. The paper \cite{Wor} proves the following result concerning the length of the longest path in $\cA_n$:
\begin{theorem}\label{th1}
There exists an absolute constant $\a$ such that if $L(n)$ denotes the length of the longest path in $\cA_n$ then
$$\Pr\brac{L(n)\geq \frac{n}{\log^\a n}}\leq \frac{1}{\log^\a n}.$$
\end{theorem}
The value of $\a$ from \cite{Wor} is rather small and we will assume for the purposes of this proof that
\beq{alpha}
\a<\frac{1}{3}.
\eeq

The aim of this paper is to give the following improvement on Theorem \ref{th1}:
\begin{theorem}\label{th2}
$$\Pr(L(n)\geq ne^{-\log^cn})\leq O(e^{-\log^{c/2}n})$$
for any constant $c<2/3$.
\end{theorem}
\proofend

This is most likely far from the truth. It is reasonable to conjecture that in fact $L(n)\leq n^{1-\e}$ w.h.p. for some positive $\e>0$. For lower bounds, \cite{Wor} shows that $L(n)\geq n^{\log_32}+2$ always and $\E(L(n))=\Omega(n^{0.8})$. Chen and Yu \cite{CY} have proved an $\Omega(n^{\log_32})$ lower bound for arbitrary 3-connected planar graphs.
\section{Outline proof strategy}
We take an arbitrary path $P$ in $\cA_n$ and bound its length. We do this as follows. We add vertices to the interior of $xyz$ in rounds.
In round $i$ we add $\s_i$ vertices. We start with $\s_0=n^{1/2}$ and choose $\s_i\gg\s_{i-1}$ where $A\gg B$ iff $B=o(A)$. We will argue inductively that $P$ only visits $\t_{i-1}=o(\s_{i-1})$ faces of $\cA_{\s_{i-1}}$ and then use Lemma \ref{P2} below to argue that roughly a fraction $\t_{i-1}/\s_{i-1}$ of the $\s_i$ new vertices go into faces visited by $P$. We then use a variant (Lemma \ref{P1}) of Theorem \ref{th1} to argue that w.h.p. $\frac{\t_i}{\s_i}\leq \frac{\t_{i-1}}{2\s_{i-1}}$. Theorem \ref{th2} will follow easily from this.

\section{Paths and Triangles}\label{pt}
Fix $1\leq \s\leq n$ and let $\cA_\s$ denote the ApN we have after inserting $\s$ vertices $A$ interior to $T_0$. It has $2\s+1$ faces, which we denote by $\cT=\set{T_1,T_2,\ldots,T_{2\s+1}}$. Now add $N$ more vertices $B$ to create a larger network $\cA_{\s'}$ where $\s'=\s+N$. Now consider a path $P=x_1,x_2,\ldots,x_m$ through $\cA_{\s'}$.
Let $I=\set{i:x_i\in A}=\set{i_1,i_2,\ldots,i_\t}$.
Note that $Q=(i_1,i_2,\ldots,i_\t)$ is a path of length $\t-1$ in $\cA_\s$. This is because $i_{k}i_{k+1}, 1\leq k<\t$ must be an edge of some face in $\cT$. We also see that for any $1\leq k<\t$ that the vertices $x_j,i_k<j< i_{k+1}$ will all be interior to the same face $T_l$ for some $l\in [2\s+1]$.

We summarise this in the following lemma: We use the notation of the preceding paragraph.
\begin{lemma}\label{lem1}
Suppose that $1\leq \s<\s'\leq n$ and that $Q$ is a path of $\cA_\s$ that is obtained from a path $P$ in $\cA_{\s'}$ by omitting the vertices in $B$.

Suppose that $Q$ has $\t$ vertices and that $P$ visits the interior of $\t'$ faces from $\cT$. Then
$$\t-1\leq \t'\leq \t+1.$$
\end{lemma}
\proofstart
The path $P$ breaks into vertices of $\cA_\s$ plus $\t+1$ intervals where in an interval it visits the interior of a single face in $\cT$. This justifies the upper bound. The lower bound comes from the fact that except for the face in which it starts, if $P$ re-enters a face $xyz$, then it cannot leave it, because it will have already visited all three vertices $x,y,z$. Thus at most two of the aforementioned intervals can represent a repeated face.
\proofend
\section{A Structural Lemma}
Let
{$$\l_1=\log^2n.$$}
\begin{lemma}\label{P2}
The following holds for all $i$. Let $\s=\s_i$ and suppose that $\l_1\leq \t\ll \s$. Suppose that $T_1,T_2,\ldots,T_\t$ is a set of triangular faces of $\cA_\s$. Suppose that $N\gg\s$ and that when adding $N$ vertices to $\cA_\s$ we find that $M_j$ vertices are placed in $T_j$ for $j=1,2,\ldots,\t$. Then for all $J\subseteq [2\s+1],\,|J|=\t$ we have
$$\sum_{j\in J}M_j\leq \frac{100\t N}{\s}\log\bfrac{\s}{\t}.$$
This holds q.s.\footnote{A sequence of events $\cE_n$ holds {\em quite surely} (q.s.) if $\Pr(\neg\cE_n)=O(n^{-K}$ for any constant $K>0$.} for all choices of $\t,\s$ and $T_1,T_2,\ldots,T_\t$.
\end{lemma}
\proofstart
We consider the following process. It is a simple example of a {\em
  branching random walk}. We consider a process that starts with $s$ newly {\em born} particles. Once a particle is born, it waits an exponentially mean one distributed amount of time. After this time, it simultaneously {\em dies} and gives birth to $k$ new particles and so on. A birth corresponds to a vertex of our network and a particle corresponds to a face.

Let $Z_t$ denote the number of deaths up to time $t$. The number of particles in the system is $\b_N=s+N(k-1)$. Then we have
$$\Pr(Z_{t+dt}=N)=\b_{N-1}\Pr(Z_t=N-1)dt+(1-\b_Ndt)\Pr(Z_t=N).$$
So, if $p_N(t)=\Pr(Z_t=N)$, we have $f_N(0)=1_{N=s}$ and
$$p_N'(t)=\b_{N-1}p_{N-1}(t)-\b_Np_N(t).$$
This yields
\begin{align*}
p_N(t)&= \prod_{i=1}^N \frac{(k-1)(i-1)+s}{(k-1)i} \; \; \times
e^{-st} (1-e^{-(k-1)t})^N\\
&=A_{k,N,s}e^{-st} (1-e^{-(k-1)t})^N.
\end{align*}
$A_{3,0,s}=1$.
When $s$ is even, $s,N\to\infty$, and $k=3$ we have
\begin{align*}
A_{3,N,s}&=\prod_{i=1}^N \bfrac{s/2+i-1}{i}=\binom{N+s/2-1}{s/2-1}\\
&\approx \brac{1+\frac{s-2}{2N}}^N\brac{1+\frac{2N}{s-2}}^{s/2-1}\sqrt{\frac{2N+s}{2\p Ns}}.
\end{align*}
We also need to have an upper bound for small even $s$, $N^2=o(s)$,
say. In this case we use
$$A_{3,N,s}\leq s^N.$$

When $s\geq 3$ is odd, $s,N\to\infty$ (no need to deal with small $N$ here) and $k=3$ we have
\begin{align*}
A_{3,N,s}&=\prod_{i=1}^N \bfrac{2i-2+s}{2i}=\frac{(s-1+2N)!((s-1)/2)!}{2^{2N}(s-1)!N!((s-1)/2+N)!}\\
&\approx   \brac{1+\frac{s-1}{2N}}^N\brac{1+\frac{2N}{s-1}}^{(s-1)/2}\frac{1}{(2\p N)^{1/2}}.
\end{align*}
We now consider with $\t\to\infty,\t\ll\s,N\geq m\geq  2\t N/\s\gg\t$ and arbitrary $t$,\\
(under the assumption that $\t$ is odd and $\s$ is odd)\\
(We sometimes use $A\leb B$ in place of $A=O(B)$).
\begin{align}
&\Pr(M_1+\cdots +M_{\t}=m\mid M_1+\cdots+M_{\s}=N)\nonumber\\
&=\frac{\Pr(M_1+\cdots +M_{\t}=m)\Pr(M_{\t+1}+\cdots +M_{\s}=N-m)}{\Pr(M_1+\cdots+M_{\s}=N)}\nonumber\\
&=\frac{A_{3,m,\t}A_{3,N-m,\s-\t}}{A_{3,N,\s}}\nonumber\\
&\approx\frac{ \brac{1+\frac{\t-1}{2m}}^m\brac{1+\frac{2m}{\t-1}}^{(\t-1)/2} \brac{1+\frac{\s-\t-2}{2(N-m)}}^{N-m}
\brac{1+\frac{2(N-m)}{\s-\t-2}}^{(\s-\t-2)/2}(N(2(N-m)+\s))^{1/2}}{\brac{1+\frac{\s-1}{2N}}^N \brac{1+\frac{2N}{\s-1}}^{(\s-1)/2}(2\p m\s(N-m))^{1/2}}\label{bound}\\
&\leb \frac{e^{(\t-1)/2}\brac{\bfrac{2m}{\t}^{(\t-1)/2}e^{o(\t)}}e^{(\s-\t)/2}\brac{1+\frac{2(N-m)}{\s-\t-2}}^{(\s-\t-2)/2}(N(2(N-m)+\s))^{1/2}}{e^{\s/2-\s^2/8N}\brac{\bfrac{2N}{\s}^{(\s-1)/2}e^{\s^2/(4+o(1))N}}(m\s(N-m))^{1/2}}\nonumber\\
&\leb \frac{e^{o(\t)}\bfrac{2m}{\t}^{(\t-1)/2}\brac{1+\frac{2(N-m)}{\s-\t-2}}^{(\s-\t-2)/2}(N(2(N-m)+\s))^{1/2}}{\bfrac{2N}{\s}^{(\s-1)/2}(m\s(N-m))^{1/2}}\nonumber
\end{align}

The above bound can be re-written as
$$\leb
\frac{e^{o(\t)}\bfrac{2}{\t}^{(\t-1)/2}N^{1/2}\s^{(\s-1)/2}}{(2N)^{(\s-1)/2}\s^{1/2}}\times \frac{m^{(\t-1)/2}\brac{1+\frac{2(N-m)}{\s-\t-2}}^{(\s-\t-2)/2}(N-m+\s)^{1/2}}{(m(N-m))^{1/2}}.$$
Suppose first that
$m\leq N-4\s$.
Then the bound becomes
\begin{align}
&\leb\frac{e^{o(\t)}\bfrac{2}{\t}^{(\t-1)/2}N^{1/2}\s^{(\s-1)/2}}{(2N)^{(\s-1)/2}\s^{1/2}}\times m^{(\t-2)/2}\brac{1+\frac{2(N-m)}{\s-\t-2}}^{(\s-\t-2)/2}\label{42}\\
&\leb\frac{e^{o(\t)}2^{(\t-1)/2}N^{1/2}\s^{(\s-1)/2}}{(2N)^{(\s-1)/2}\t^{\t/2}}\times m^{(\t-2)/2} \bfrac{2(N-m)}{\s-\t}^{(\s-\t)/2}e^{\s^2/(N-m)}\nonumber\\
&\leq \frac{e^{o(\t)}N^{1/2}}{m^{1/2}}\bfrac{\s(N-m)}{N(\s-\t)}^{(\s-\t)/2}\bfrac{\s m}{\t N}^{(\t-1)/2}e^{\s^2/(N-m)}\nonumber\\
&\leb \frac{e^{o(\t)}N^{1/2}}{m^{1/2}}\brac{\frac{e^2m\s}{\t N}\cdot
\exp\set{-\frac{m(\s-\t)}{(\t-1)N}+\frac{2\s^2}{(\t-1)(N-m)}}}^{(\t-1)/2}.\nonumber\\
&=\frac{e^{o(\t)}N^{1/2}}{m^{1/2}}\brac{\frac{e^2m\s}{\t N}
\cdot \exp\set{-\frac{m\s}{(\t-1)N}
\brac{1-\frac{\t}{\s}-\frac{2\s}{m}-\frac{2\s}{N-m}}}}^{(\t-1)/2}\nonumber\\
&\leq \frac{e^{o(\t)}N^{1/2}}{m^{1/2}}\brac{\frac{e^2m\s}{\t N}
\cdot \exp\set{-\frac{m\s}{3\t N}}}^{(\t-1)/2}\nonumber
\end{align}
We inflate this by $n^2\binom{2\s+1}{\t}$ to account for our choices for $\s,\t,T_1,\ldots,T_\t$ to get
$$\leb n^2\frac{e^{o(\t)}N^{1/2}}{m^{1/2}}\brac{\frac{4e^4m\s^3}{\t^3 N}
\cdot \exp\set{-\frac{m\s}{3\t N}}}^{(\t-1)/2}.$$
So, if $m_0=\frac{100\t N\log(\s/\t)}{\s}$ then
\begin{align*}
&\sum_{m=m_0}^{N-4\s}\Pr(\exists \s,\t,T_1,\ldots,T_\t:M_1+\cdots +M_{\t}=m\mid M_1+\cdots+M_{\s}=N)\\
&\leb n^2e^{o(\t)}N^{5/2}\sum_{m=m_0}^{N-4\s}\brac{\frac{4e^4m\s^3}{\t^3 N}\cdot \exp\set{-\frac{m\s}{3\t N}}}^{(\t-1)/2}\\
&\leq n^2e^{o(\t)}N^{7/2}\brac{\frac{4e^4m_0\s^3}{\t^3 N}\cdot \exp\set{-\frac{m_0\s}{3\t N}}}^{(\t-1)/2}\\
\noalign{since $xe^{-Ax}$ is decreasing for $Ax\geq 1$}\\
&=n^2e^{o(\t)}N^{7/2}\brac{\frac{4e^4m_0\s}{\t N}\exp\set{-\frac{m_0\s}{6\t N}}\times \frac{\s^2}{\t^2}\exp\set{-\frac{m_0\s}{6\t N}}}^{(\t-1)/2}\\
&\leq n^2N^{7/2}\brac{400e^{4+o(1)}\log\bfrac{\s}{\t}\times e^{-50/3}\times \frac{\s^2}{\t^2}\bfrac{\t}{\s}^{50/3}}^{(\t-1)/2}\\
&=O(n^{-anyconstant}).
\end{align*}
Suppose now that $N-4\s\leq m\leq N-\s^{1/3}$. Then we can bound \eqref{42} by
\begin{align*}
&\leb \frac{e^{o(\t)}\bfrac{2}{\t}^{(\t-1)/2}\s^{(\s-1)/2}}{(2N)^{(\s-1)/2}}\times m^{(\t-1)/2}e^{4\s}\\
&\leq \bfrac{e^8\s}{2N}^{(\s-\t)/2}\bfrac{e^8\s}{\t}^{(\t-1)/2}.
\end{align*}
We inflate this by $n^2\binom{2\s+1}{\t}<n^24^\s$ to get
$$\leb n^2\bfrac{8e^8\s}{N}^{(\s-\t)/2}\bfrac{16e^8\s}{\t}^{(\t-1)/2}$$
So,
\begin{align*}
&\sum_{m=N-4\s}^{N-\s^{1/3}}\Pr(\exists \s,\t,T_1,\ldots,T_\s:M_1+\cdots +M_{\t}=m\mid M_1+\cdots+M_{\s}=N)\\
&\leb n^2N^2\s \bfrac{8e^8\s}{N}^{(\s-\t)/2}\bfrac{16e^8\s}{\t}^{(\t-1)/2}\\
&=O(n^{-anyconstant})
\end{align*}
since $\s\log N\gg\t\log\s$.

When $m\geq N-\s^{1/3}$ we replace \eqref{bound} by
\begin{align*}
&\leb\frac{ \brac{1+\frac{\t-1}{2m}}^m\brac{1+\frac{2m}{\t-1}}^{(\t-1)/2}\s^{N-m}
N^{1/2}}{\brac{1+\frac{\s-1}{2N}}^N
\brac{1+\frac{2N}{\s-1}}^{(\s-1)/2}(m\s)^{1/2}}\\
&\leb
\frac{e^{\t/2+o(\t)}\bfrac{2m}{\t}^{(\t-1)/2}\s^{N-m}N^{1/2}}{e^\s\bfrac{2N}{\s}^{(\s-1)/2}m^{1/2}}\\
&\leb
\bfrac{e^{1+o(1)}\s}{\t}^{(\t-1)/2}\bfrac{\s}{2N}^{(\s-\t)/2}\s^{\s^{1/3}}.
\end{align*}
Inflating this by $n^24^\s$ gives a bound of
$$\leb n^2\bfrac{16e^{1+o(1)}\s}{\t}^{(\t-1)/2}\bfrac{8\s^{1+o(1)}}{N}^{(\s-\t)/2}= O(n^{-anyconstant}).$$

\proofend

\section{Modifications of  Theorem \ref{th1}}\label{mod}
Let $\l=\log^3n$ and partition $[\l]$ into $q=\log n$ sets of size $\l_1=\log^2 n$. Now add $n-\l$ vertices to $\cT_\l$ and let $M_i$ denote the number of vertices that land in the $i$th part $\Pi_i$ of the partition. Lemma \ref{P2} implies that q.s.
\beq{FUKU}
M_i\leq { M_{\max}=}\frac{200n}{\log n}\log\log n,\quad 1\leq i\leq \t.
\eeq
Let
\beq{omega}
\om_1(x)=\log^{\a/2}x
\eeq
for $x\in \Re$.

Let $L_i$ denote the length of the longest path in $\Pi_i$. Suppose that $\cT_n$ contains a path of length at least $n/\om_1,\,\om_1=\om_1(n)$ and let $k$ be the number of $i$ such that
{$$L_i\geq \frac{200n\log\log n}{\om_1^2\log n}\geq \frac{M_{\max}}{\log^\a(M_{\max})}.$$}
Then, as $k \le q=\log n$ we have
$$k\frac{200n\log\log n}{\log n} +(\log n-k)\frac{200n\log\log n}{\om_1^2\log n}\geq \frac{n}{\om_1}$$
which implies that 
$$k\geq \frac{\log n}{201\om_1\log\log n}.$$
Theorem \ref{th1}
with the bound on $M_i$ given in \eqref{FUKU}
implies that the probability of this is at most
\beq{fwit}
\frac{1}{n}+\binom{\log n}{\frac{\log n}{201\om_1\log\log n}}\bfrac{1}{\log^{\a}(n/\log n)}^{\frac{\log n}{201\om_1\log\log n}}\leq\frac{1}{n}+
\bfrac{1}{\log^{\a/3}n}^{\frac{\log n}{201\om_1\log\log n}}
\leq \frac1{\f(n,\om_1)}
\eeq
where
$$\f(x,y)=\exp\set{\frac{\log x}{y\log\log x}}.$$
The term $1/n$ accounts for the failure of the property in Lemma \ref{P2}.

In summary, we have proved the following
\begin{lemma}\label{P1}
\beq{P3}
\Pr\brac{L(n)\geq \frac{n}{\om_1(n)}}\leq \frac{1}{\f(n,\om_1)}.
\eeq
\end{lemma}
\proofend

We are using $\f(x,y)$ in place of $\f(x)$ because we will need to use $\om_1(x)$ for values of $x$ other than $n$.

Next consider $\cA_\s$ and $\l_1\leq \t\ll\s$ and let $T_1,T_2,\ldots,T_\t$ be a set of $\t$ triangular faces of $\cA_\s$. Suppose that we add $N\gg\s$ more vertices and let $N_j$ be the number of vertices that are placed in $T_j,\,1\leq j\leq \t$.

Next let
\beq{LL}
\LL(x)=e^{x^2}
\eeq
where $x\in \Re$.

Now let
\beq{SeeJ}
J=\set{j:N_j\geq \LL_0}\text{ where }\LL_0=\LL(\om_1(n)).
\eeq
Let $L_j$ denote the length of the longest path through the ApN
defined by $T_j$ and the $N_j$ vertices it contains, $1\leq j\leq
\t$. {For the remainder of the section let
\beq{om2}
\om_0=\om_1(\LL_0),\quad\f_0=\f(\LL_0,\om_0)=\exp\set{\frac{\om_0}{2\log\om_0}},\quad\om_2=\frac{\f_0}{\om_0}.
\eeq}
Then let
\beq{J1}
J_1=\set{j\in J:L_j\geq \frac{N_j}{\om_1(N_j)}}.
\eeq

We note that
\begin{align*}
\log\om_2&=\log\f_0-\log\om_0=\frac{\log\LL_0}{\om_0\log\log\LL}-\log\om_0\\
&=\frac{\om_0^2}{(2+o(1))\om_0\log\log\om_0}-\log\om_0.
\end{align*}
For $j \in J$, $N_j  \ge \LL_0$ (see \eqref{SeeJ}).
It  follows from Lemma \ref{P1} that the size of $J_1$ is stochastically
dominated by $Bin(\t,1/\f_0)$. Using a  Chernoff bound we  find that
{\beq{P8}
\Pr\brac{|J_1|\geq \frac{\om_2\t}{\f_0}}\leq \bfrac{e}{\om_2}^{\om_2\t/\f_0}.
\eeq}

Using this we prove
\begin{lemma}\label{P7}
Suppose that
$$\log\bfrac{\s}{\t}\leq \frac{\om_0}{\log\om_0}.$$
Then q.s., for all $\l_1\leq\t\ll\s\ll N$ and all collections $\cT$ of $\t$ faces of $\cA_\s$ we find that with $J_1$ as defined in \eqref{J1},
{$$|J_1|\leq \frac{\om_2\t}{\f_0}.$$}
\end{lemma}
\proofstart
It follows from \eqref{P8} that
\begin{align*}
&\Pr\brac{\exists \t,\s,N,\cT:|J_1|\geq\frac{\om_2}{\t\f_0}} \\
&\leq n^3\binom{(2\s+1)}{\t}\bfrac{e}{\om_2}^{\om_2\t/\f_0}\\
&\leq n^3\brac{\frac{e(2\s+1)}{\t}\cdot\bfrac{e}{\om_2}^{\om_2/\f_0}}^\t\\
&\leq \exp\set{\t\brac{\frac{3\log n}{\t}+2+\log\bfrac{\s}{\t}+
\frac{\om_2}{\f_0}-\frac{\om_2\log\om_2}{\f_0}}}\\
&\leq \exp\set{\t\brac{\frac{3\log n}{\t}+2+\frac{\om_0}{\log\om_0}+ {-\frac{\om_0}{(2+o(1))\log\log\om_0}}}}\\
&=O(n^{-anyconstant}).
\end{align*}
\proofend
\section{Proof of Theorem \ref{th2}}\label{pot2}
Fix a path $P$ of $\cA_n$. Suppose that after adding
$\s\geq n^{1/2}$ vertices we find that $P$ visits
{\beq{t1f}
n^{1/2}\geq \t\geq \l_1\om_0
\eeq}
of the triangles $T_1,T_2,\ldots,T_\t$ of $\cA_\s$.
Now consider adding $N$ more vertices, where the value of $N$ is given in \eqref{N=} below.
Let $\s'=\s+N$ and let $\t'$ be the number of triangles of $\cA_{\s'}$ that are visited by $P$.

We assume that
\beq{P5}
\frac{\a}{2}\log\log n\leq \log\bfrac{\s}{\t}\leq \frac{\om_0}{\log\om_0}.
\eeq
Let $M_i$ be the number of vertices placed in $T_i$ and let $N_i$ be the number of these that are visited by $P$. It follows from Lemma \ref{P2} that w,h.p.
$$\sum_{i=1}^\t M_i\leq \frac{100\t N}{\s}\log\bfrac{\s}{\t}.$$

Now w.h.p.,
\beq{P6}
\sum_{i=1}^\t N_i\leq \t\LL_0+\frac{100\om_2\t N}{\f_0\s}\log\bfrac{\s\f_0}{\om_2\t}+
\frac{100\t N}{\s\om_0}\log\bfrac{\s}{\t}.
\eeq
{\bf Explanation:} $\t\LL_0$ bounds the contribution from $[\t]\setminus J$ (see \eqref{SeeJ}). The second term bounds the contribution from $J_1$.
{Now $|J_1|<\om_2\t/\f_0\ll\t$ as shown in Lemma \ref{P7}. We
  cannot apply}
{Lemma \ref{P2} to bound the contribution of $J_1$ unless we know
  that $|J_1|\geq \l_1$. We choose an arbitrary set of indices $J_2\subseteq
  [\t]\setminus J_1$ of size $\om_2\t/\f_0-|J_1|$ and then the middle
  term bounds the contribution of $J_1\cup J_2$. Note that $\om_2\t/\f_0=\t/\om_0\geq
\l_1$ from \eqref{t1f}.} The third term bounds the contribution from
$J\setminus J_1$.
{Here we use $\om_1(N_j)\geq \om_1(\LL_0)=\om_0$, see \eqref{J1}.}

We now  choose
\beq{N=}
N=3\s\LL_0.
\eeq
We observe that
$$\frac{\om_2}{\f_0}\log\bfrac{\s\f_0}{\om_2\t}\leq \frac{1}{\om_0}\brac{\frac{\om_0}{\log\om_0}+2\log\om_0}=o(1).$$
$$\frac{1}{\om_0}\log\bfrac{\s}{\t}\leq\frac{1}{\log\om_0}=o(1).$$
Now along with Lemma \ref{lem1} this implies that
$$\t'\leq \sum_{i=1}^\t(N_i+1)\leq \t+\t\LL_0+o\bfrac{\t N}{\s}.$$
Since $\s'=\s+N$ this implies that
$$\frac{\t'}{\s'}\leq
\brac{\frac{1}{3}+o(1)}\frac{\t}{\s}<\frac{\t}{2\s}.$$
It follows by repeated application of this argument that we can
replace Theorem \ref{th1} by
\begin{lemma}\label{P4}
{$$\Pr\brac{L(n)\geq \log n+\frac{100\log n}{e^{\om_0/\log\om_0}}n}=O\bfrac{1}{\f(n,\om_1(n))}.$$}
\end{lemma}
\proofstart
{We add the vertices in rounds of size $\s_0=n^{1/2},\s_1,\ldots,\s_m$. Here
$\s_i=3\s_{i-1}\LL_0$ and $m-1\geq (1-o(1))\frac{\log
  n}{\log\LL_0}=(1-o(1))\frac{\log n}{\om_1(n)^2}=\log^{1-2\a}n$. We let $P_0,P_1,P_2,\ldots,P_m=P$ be a sequence of paths where $P_i$ is a path in $\cA_i=\cA_{\s_0+\cdots+\s_i}$. Furthermore, $P_i$ is obtained from $P_{i+1}$ in the same way that $Q$ is obtained from $P$ in Lemma \ref{lem1}. We let $\t_i$} 
denote the number of faces of $\cA_i$ whose interior is visited by
  $P_i$. It follows from Lemma \ref{lem1} and Lemma \ref{P2} that the
  length of $P$ is bounded by 
$$m+\frac{\t_{m-1}}{\s_{m-1}}\s_m\log\bfrac{\s_{m-1}}{\t_{m-1}},$$
since the second term is a bound on the number of points in the
interior of triangles of $\cA_{m-1}$ visited by $P$.

We have w.h.p. that
$$\frac{\s_i}{\t_i}\geq \begin{cases}\frac{2\s_{i-1}}{\t_{i-1}}&\frac{\s_{i-1}}{\t_{i-1}}\leq
  e^{\om_0/\log\om_0}\\ \frac{\s_{i-1}}{100\t_{i-1}\log(\s_{i-1}/\t_{i-1})}&\frac{\s_{i-1}}{\t_{i-1}}>
  e^{\om_0/\log\om_0}\end{cases}.$$
The second inequality here is from Lemma \ref{P2}.

The result follows from $2^{\log^{1-2\a}n}\geq e^{\om_0/\log\om_0}$.

\proofend

To get Theorem \ref{th2} we repeat the argument in Sections \ref{mod} and
\ref{pot2}, but we start with $\om_1(x)=\log^{1/3}x$. The claim in
Theorem \ref{th2} is then slightly weaker than the claim in Lemma \ref{P4}.

\end{document}